\input amstex
\input amsppt.sty

\def\epsilon{\varepsilon}
\def\capa{\operatorname{cap}}

\def\phi{\varphi}

\magnification=1200

\NoBlackBoxes
        \topmatter
\title
Bergman completeness is not a quasi-conformal invariant
\endtitle
\author
Xu Wang
\endauthor
\abstract
We show that Bergman
completeness is not a quasi-conformal invariant for general Riemann surfaces.
\endabstract
\address
Instytut Matematyki, Uniwersytet Jagiello\'nski, \L ojasiewicza 6, 30-348
Krak\'ow, Po\-land
\endaddress
\email
xu.wang\@im.uj.edu.pl
\endemail
\thanks
I would like to thank W. Zwonek for his fruitful suggestions on this paper. Project operated within the Foundation for Polish Science
IPP Programme "Geometry and Topology in Physical Models"
co-financed by the EU European Regional Development Fund,
Operational Program Innovative Economy 2007-2013
\endthanks

\endtopmatter

\document

Recently, B.-Y. Chen (see \cite{Chen~1}) asked a question whether Bergman
completeness is a quasi-conformal invariant of Riemann surfaces. In this paper we will
give an example showing that the answer to the question is negative.

\medskip

The original idea of the example is based on the paper \cite{Zwo~1}. In \cite{Zwo~1}, the growth of the Bergman kernel
near the boundary has been estimated with the help of potential-theoretical quantities. Also in \cite{Ju},
a necessary and sufficient condition for Bergman completeness of Zalcman type domains has been found.
The present paper will use the
methods in \cite{Pfl-Zwo}.

\medskip

For $0<r<\frac{1}{4}$ and $t\in\left (0,\frac{1}{2}\right)$ define
$$
A^{r,t}:=\bigcup_{k=1}^{\infty}A^{r,t}_{k}\cup \{0\}
$$
where
$$
A^{r,t}_{k}:=\{r^{k}e^{i\theta}: -2\alpha_{k} \leq \theta \leq 2\alpha_{k}\}
$$
and $\sin \alpha_{k}=e^{-t^{-k}}$.

Finally we put $D^{r,t}:=\triangle(0,1)\setminus A^{r,t}$ where $\triangle(p,r):=\{z\in\Bbb C:|z-p|<r\}$, $p\in\Bbb C$, $r>0$.

\proclaim{Theorem 1}
\medskip
\item{(1)} $D^{r,t}$ is Bergman exhaustive if and only if $r^{2}\leq t$;
\item{(2)} if $t\leq r^{4}$ then $D^{r,t}$ is not Bergman complete.
\endproclaim

\medskip

\proclaim{Corollary 2}
Bergman completeness is not in general a quasi-conformal invariant for Riemann surfaces.
\endproclaim

\medskip

Recall that a homeomorphism $f$ defined on a domain in the complex plane is called $L$ ($L>1$) quasi-conformal
if it is differentiable almost everywhere and
$$
\left |\frac{\partial f}{\partial \overline{z}}\right| \leq \frac{L-1}{L+1}\left|\frac{\partial f}{\partial z}\right|.
$$

Before we present the proofs of Theorem 1 and Corollary 2 let us recall some notions and results which we need in the paper.
\medskip

In order to prove Theorem 1, we need some lemmas on logarithmic capacity and Bergman kernel (see \cite{Jar-Pfl}). First we recall necessary notions
from potential theory (see \cite{Ran} where also other properties of potential-theoretic objects that we use are given).

\medskip

For a probabilistic measure $\mu$ defined on all Borel subsets of a compact set
$K$ (denote $\mu\in\Cal P(K)$), we define its logarithmic potential
by
$$
p_{\mu}(z):=\int_K\log|z-w|d\mu(w),\;z\in\Bbb C.
$$
Recall that $p_{\mu}$ is harmonic in $\Bbb C\setminus K$ and subharmonic in $\Bbb C$.

\medskip

We also denote the energy of $\mu$ as follows
$$
I(\mu):=\int_Kp_{\mu}(z)d\mu(z)=\int_K\int_K\log|w-z|d\mu(w)d\mu(z).
$$

\medskip

A probabilistic measure $\nu$ defined on Borel subsets of the
compact set $K$ is called {\it the equilibrium measure of $K$} if
$I(\nu)=\sup\{I(\mu):\mu\in\Cal P(K)\}$. It is well-known that the
equilibrium measure exists and is unique if $K$ is not polar (a
set $F\subset \Bbb C$ is called polar if there is a subharmonic
function $u$ defined on $\Bbb C$ such that $u\not\equiv-\infty$
and $F\subset\{u=-\infty\}$).

\medskip

{\it The logarithmic capacity} of a subset $E$ of $\Bbb C$ is
given by the formula
$$
\capa(E):=e^{\sup\{I(\mu):\mu\in\Cal P(K), K\text{is a compact subset of $E$}\}}.
$$

\medskip

In case when $K$ is
compact and not polar then $\capa(K)=e^{I(\nu)}$, where $\nu$
denotes the equilibrium measure of $K$. It is well-known that a
Borel set $E\subset\Bbb C$ is polar iff $\capa(E)=0$.

\bigskip

Now let us recall the notion of the Bergman kernel. Let
$D$ be a bounded domain in $\Bbb C^n$. Denote by $L_h^2(D)$
square integrable holomorphic functions on $D$. $L_h^2(D)$ is a
Hilbert space with the scalar product induced from $L^2(D)$
denoted by $\langle\cdot,\cdot\rangle_D$ (we also denote by $||\cdot||_D$ the $L^2$ norm on $D$ and the space of square integrable holomorphic functions on $D$ is denoted by $L_h^2(D)$).

Let us define the {\it Bergman kernel of $D$} as
$$
K_D(z)=\sup\{\frac{|f(z)|^2}{||f||_D^2}:f\in L^2_h(D), f\not\equiv 0\}, z\in D
$$
and the fundamental form of the Bergman metric as
$$
B_D(z)=i\partial\overline{\partial}\log K_D(z), z\in D
$$

\medskip

Let
$$
I_D(z,X)=\sup\{\frac{|f^{\prime}(z)X|^2}{||f||_D^2}:f\in L^2_h(D), f(z)=0, f\not\equiv 0\}, z\in D, X\in \Bbb C^n.
$$
The following result is classical
$$
B_D(z)(X)=\frac{I_D(z,X)}{K_D(z)}, z\in D, X\in \Bbb C^n.
$$

Let $D$ be a bounded domain in $\Bbb C^n$, $z_0\in\partial D$.
Then $D$ is {\it Bergman exhaustive  at $z_0\in\partial D$ } if
$\lim_{D\owns z\to z_0}K_D(z)=\infty$. We call $D$  {\it Bergman
exhaustive } if $D$ is Bergman exhaustive at $z_0$ for any
$z_0\in\partial D$.

\medskip

A bounded domain $D$ is called {\it Bergman complete} if any
Cauchy sequence with respect to the Bergman distance is convergent to some point in $D$
under the standard topology of $D$. For references on the
Bergman kernel, metric and distance, see \cite{Chen~2}.

\medskip

It is known that if a bounded domain
$D\subset\Bbb C$ is Bergman exhaustive then $D$ is Bergman
complete (see \cite{Chen~3}), the converse implication does not
hold in general (see \cite{Zwo~2}). Also bounded hyperconvex domain
is Bergman exhaustive (see \cite{Ohs})
and Bergman complete (see \cite{B-P}, \cite{Her}).
\bigskip

Now let us introduce the notions necessary for the description of
Bergman exhaustive points in dimension one.

\medskip

For a bounded domain $D\subset\Bbb C$ and for a point $z\in\bar D$, we introduce the following potential theoretic quantity;
$$
\gamma_D^{(n)}(z):=\int_0^{1/4}\frac{d\delta}{\delta^{2n+3}(-\log(\capa(\bar\triangle(z,\delta)\setminus D))}.
$$

\medskip

The following lemma comes from the paper \cite{Zwo~1}.

\proclaim{Lemma 3} Let $D$ be a bounded domain in $\Bbb C$ and let $z_0\in\partial D$. Then
$$
\lim_{D\owns z\to z_0}\gamma_D^{(0)}(z)=\infty        \tag{1}
$$
if and only if
$$
\text{ $D$ is Bergman exhaustive at $z_0$.}      \tag{2}
$$
\endproclaim

\subheading{Remark 4}
It also follows from classical results that the domain $D^{r,t}$ is Bergman exhaustive at all of its boundary point except for $0$.

\bigskip

Since we shall only consider bounded domains in $\Bbb C$, for simplicity we denote $\beta_{D}(z)=B_{D}(z)(1),z\in D$. From the paper \cite{Pfl-Zwo}, we know

\proclaim{Lemma 5} Let $D$ be a bounded domain in $\Bbb C$, $D\owns z_{k}\to z_0\in\partial D$. If
$$
\limsup\sb{k\to\infty}\gamma_D^{(1)}(z_{k})<\infty,
$$
then
$$
\limsup\sb{k\to\infty}\beta_D(z_{k})<\infty.
$$
\endproclaim

Now we move to the proofs of main results. Let us first see how we derive Corollary 2 from Theorem 1. Then we show Theorem 1.

\bigskip
\demo{Proof of Corollary 2}
For $\alpha>\frac{1}{2}$ we define
$$
\varphi(z)=z^{\alpha}\overline{z}^{\alpha-1}.
$$
Note that $\varphi$ is a quasi-conformal mapping from $D^{r,t}$ to $D^{r^{2\alpha-1},t}$.
Choosing for instance $r=1/8$, $\alpha=2/3$, $t=1/32$ we get that Bergman
completeness is not a quasi-conformal invariant. Actually $D^{r,t}$ is, in view of Theorem 1 (1),  Bergman exhaustive and thus it is Bergman complete whereas $D^{r^{2\alpha-1},t}$ is, in view of Theorem 1 (2), not Bergman complete.
\qed
\enddemo

\subheading{Remark 6}
By Theorem 1 we also can deduce that Bergman exhaustiveness is not a quasi-conformal invariant.
But we still do not know whether it is a conformal invariant for bounded domains in the complex plane.
Also due to the above theorem, there are lots of domains which are Bergman complete but not
Bergman exhaustive (see \cite{Zwo~2}).

Note also that the example from Corollary 2 is a Riemann surface with infinite dimensional
fundamental group. But B.-Y. Chen informed the author that Bergman
completeness is a quasi-conformal invariant when the fundamental group is finitely generated.

Behind Chen's question on the quasi-conformal
invariance of the Bergman completeness, there is an affirmative
result by Pfluger on the quasi-conformal invariance of the existence
of the Green function, and a very deep
negative result by Beurling-Ahlfors on the invariance of the nullity
of linear measure for subsets of the circle as the boundary of the
disc (see \cite{Sa-Nak}). The author would like to thank the referee for this comment.

\demo{Proof of Theorem 1}

If $D^{r,t}$ is Bergman exhaustive then by Lemma 3,
$$
\lim_{D^{r,t}\owns z\to 0}\gamma_{D^{r,t}}^{(0)}(z)=\infty.
$$
For any $\delta\in [0,1/4]$, if $-1<x<0$ and $|x|$ small enough we have
$$
\bar\triangle(x,\delta)\setminus D^{r,t} \subset \bar\triangle(0,\delta)\setminus D^{r,t}.
$$
Thus for such $x$
$$
\gamma_{D^{r,t}}^{(n)}(x)\leq \gamma_{D^{r,t}}^{(n)}(0),
$$
consequently, we get $\gamma_{D^{r,t}}^{(0)}(0)=\infty$.
Now choose
$$
A_{1}=\{z\in \Bbb C; 2r^{2} \leq |z| \leq 1/4 \},
$$
$$
A_{k}=\{z\in \Bbb C; 2r^{k+1} \leq |z| \leq 2r^{k} \}, k \geq 2.
$$
We have
$$
\gamma_{D^{r,t}}^{(n)}(0)=(\int_{2r^{2}}^{1/4}+\sum_{k=2}^{\infty}\int_{2r^{k+1}}^{2r^{k}})
\frac{d\delta}{\delta^{2n+3}(-\log(\capa(\bar\triangle(0,\delta)\setminus D^{r,t}))}.
$$
Let
$$
C_{1}=\int_{2r^{2}}^{1/4}\frac{d\delta}{\delta^{2n+3}(-\log(\capa(\bar\triangle(0,\delta)\setminus D^{r,t}))},
$$
$$
C_{k}=\int_{2r^{k+1}}^{2r^{k}}\frac{d\delta}{\delta^{2n+3}(-\log(\capa(\bar\triangle(0,\delta)\setminus D^{r,t}))}, k \geq 2.
$$
Then we have
$$
(1/4-2r^{2})4^{2n+3}\frac{1}{-\log(\capa(A_{2}\setminus D^{r,t}))} \leq C_{1},
$$
$$
C_{1}\leq (1/4-2r^{2})(2r^{2})^{-2n-3}\sum_{j=1}^{\infty}\frac{1}{-\log(\capa(A_{j}\setminus D^{r,t}))},
$$
$$
(2r^{k})^{-2n-2}(1-r)\frac{1}{-\log(\capa(A_{k+1}\setminus D^{r,t}))}\leq C_{k},
$$
$$
C_{k}\leq (2r^{k+1})^{-2n-2}(1/r-1)\sum_{j=k}^{\infty}\frac{1}{-\log(\capa(A_{j}\setminus D^{r,t}))},k\geq 2.
$$
As $\capa(A_{j}\setminus D^{r,t})=r^{j}e^{-t^{-j}}$, we get
$$
\frac{t^{j}}{1-t\log r}\leq\frac{1}{-\log(\capa(A_{j}\setminus D^{r,t}))}\leq t^{j}.
$$

\medskip

Thus there is a sufficiently large constant $C(t,r,n)$ such that

\medskip

$$
C(t,r,n)^{-1}(\frac{t}{r^{2n+2}})^{k} \leq C_{k} \leq C(t,r,n)(\frac{t}{r^{2n+2}})^{k}.
$$

\medskip
Now due to the lower semi-continuity of $\gamma_{D^{r}}^{(0)}$ (\cite{Zwo~1}), the theorem follows easily from Lemma 3 (part (1)) and Lemma 6 (part (2)).\qed\enddemo

\enddocument